\documentclass[11pt]{amsart}

\usepackage{a4}

\usepackage{amssymb,latexsym}

\usepackage[all]{xy}

\usepackage[usenames]{color}

\usepackage[normalem]{ulem}

\usepackage{hyperref}

\usepackage{mathrsfs}

\usepackage{amssymb,latexsym}

\usepackage[all]{xy}

\usepackage{enumerate}

\newcommand{\N}{\mathbb{N}}
\newcommand{\Z}{\mathbb{Z}}

\newcommand{\CC}{\mathscr{C}}

\newcommand{\CM}{\mathscr{M}}
\newcommand{\CP}{\mathscr{P}}

\newcommand{\CT}{\mathcal{T}}

\newcommand{\tasu}{T_{(A,\s,u)}}
\newcommand{\sas}{\Sym(A,\s)}
\newcommand{\das}{D_{(A,\s)}}

\newcommand{\ox}{\otimes}
\newcommand{\x}{\times}

\newcommand{\vf}{\varphi}
\newcommand{\ve}{\varepsilon}
\newcommand{\vt}{\vartheta}
\newcommand{\s}{\sigma}

\newcommand{\simtoo}{\overset\sim\longrightarrow}

\newcommand{\df}{\emph}

\newcommand{\id}{\mathrm{id}}

\DeclareMathOperator{\im}{Im}

\DeclareMathOperator{\sign}{sign}

\DeclareMathOperator{\Sym}{Sym}

\DeclareMathOperator{\Hom}{Hom}
\DeclareMathOperator{\Nil}{Nil}
 
\DeclareMathOperator{\coker}{coker}
\DeclareMathOperator{\Trd}{Trd}

\renewcommand{\leq}{\leqslant}

\newcommand{\pf}[1]{\langle\!\langle #1\rangle\!\rangle}

\newcommand{\qf}[1]{\langle #1\rangle}

\newcommand{\sm}{\setminus}
\newcommand{\wt}{\widetilde}

\DeclareMathOperator{\PSD}{PSD}
\DeclareMathOperator{\Tr}{\mathrm{Tr}}
\DeclareMathOperator{\tr}{tr}

\newtheorem{lemma}{Lemma}[section]
\newtheorem{thm}[lemma]{Theorem}
\newtheorem{prop}[lemma]{Proposition}
\newtheorem{cor}[lemma]{Corollary}

\theoremstyle{definition}
\newtheorem{defi}[lemma]{Definition}

\newtheorem{remark}[lemma]{Remark}

\newtheorem{examples}[lemma]{Examples}

\numberwithin{equation}{section}

\date{April 9, 2018}

\title[Signatures, sums of hermitian squares and  positive cones]{Signatures, sums of hermitian squares  
and positive cones on algebras with involution}
\author{Vincent Astier}
\author{Thomas Unger}
\address{School of Mathematics and Statistics\\ University College Dublin\\ Belfield\\
Dublin~4\\ Ireland} 
\email{vincent.astier@ucd.ie, thomas.unger@ucd.ie}



\begin{document}

\maketitle

\begin{abstract}
We provide a coherent picture of our  efforts thus far in extending real algebra and its links to
the theory of quadratic forms 
over ordered fields in the noncommutative direction, using hermitian forms and ``ordered'' algebras with involution.
\end{abstract}

\section{Introduction}

Hilbert, in his famous Paris talk in 1900 \cite{Hilbert1900} posed the following as his 17th problem: \emph{ob nicht jede 
definite Form als Quotient von Summen von Formenquadraten dargestellt werden kann}. Translated: is every 
nonnegative $n$-ary polynomial  over a field $F$ a sum of squares of rational functions over $F$? Partial answers were known to Hilbert, as well as the fact that a positive semidefinite polynomial need not be a sum of squares of polynomials.

The full, affirmative, answer was obtained by Artin in 1927, and built on joint work with Schreier, published in the same 
year. To be precise: let $F$ be a field of characteristic different from $2$ with space of orderings $X_F$. 
The Artin-Schreier theorem \cite{Artin-Schreier} says
that   $F$ admits an ordering 
if and only if $-1$ is not a sum of squares in $F$.  Artin's theorem \cite{Artin} says that
an element is positive at all $P\in X_F$ if and only
if it is a sum of squares in $F$. This result is crucial in his solution of Hilbert's 17th problem.

Quadratic forms over $F$ come into the picture via the Witt ring $W(F)$ (the ring of Witt equivalence
classes of nondegenerate quadratic forms over $F$). Orderings $P\in X_F$ correspond 
to signature homomorphisms $\sign_P: W(F)\to \Z$. This goes back to Sylvester's \emph{Law of Inertia} \cite{Sylvester1852}.
In addition to the fundamental ideal of $W(F)$, the prime ideals of $W(F)$ are 
given by $\sign_P^{-1}(0)$ and $\sign_P^{-1}(p\Z)$ (where $p$ is an odd prime) with $P\in X_F$ (Lorenz-Leicht 
\cite{Lorenz-Leicht}). 
Given $q\in W(F)$, the total signature map
$\sign_\bullet q: X_P \to \Z, P\mapsto \sign_P q$ is continuous with respect to the Harrison topology on $X_F$ and the discrete 
topology on $\Z$. 
Furthermore, the torsion part of $W(F)$ consists of those $q$ with $\sign q=0$ (Pfister's local-global principle 
\cite{Pfister66}).

The above results (well-documented in \cite{Lam} and \cite{Sch}), and some of their extensions to commutative rings, are among the foundations of real algebra, see for example \cite{BCR} or
Lam's expository paper \cite{Lam-1984}. In a series of recent papers \cite{A-U-Kneb, A-U-prime, A-U-pos,
 A-U-PS, A-U-stab}
(and also  \cite{LU1}) we extended these results in the noncommutative direction, more precisely
to central simple $F$-algebras with involution and hermitian forms over such algebras.

The study of central simple
algebras with involution was initiated  by Albert in the 1930s \cite{Albert35} and is still a topic of current research  as testified by \emph{The Book of Involutions}
\cite{BOI}; see also  \cite{auel-2011} and the copious references therein  for a list of open problems in this area. A large part of present day 
research in algebras with involution is driven by the deep connections with linear algebraic groups, first observed by Weil 
\cite{Weil}; see also Tignol's \emph{2 ECM} exposition \cite{Tignol-1998}. Some work has been done on algebras with involution over formally 
real fields, for example \cite{LT}, \cite{Q}, but this part of the theory is relatively underdeveloped. This observation, together with the fact that
algebras with involution are a natural generalization of quadratic forms, are motivating factors for our research. 

This article is an expanded version of the prepublication \cite{A-U-survey2017}, from the
\emph{S\'emi\-naire de Structures Alg\'ebriques Ordonn\'ees}, Universities Paris 6 and 7.

\section{Signatures}\label{sec:two}

Let $(A,\s)$ be an \emph{$F$-algebra with involution}, by which we mean  that $A$ is a finite dimensional simple $F$-algebra
with centre a field $K\supseteq F$  and $\s$ is an $F$-linear anti-automorphism of $A$ of order $2$ (which
implies that $[K:F]\leq 2$).
Let $W(A,\s)$ denote the Witt group of $(A,\s)$, i.e. the $W(F)$-module of Witt equivalence classes of 
nondegenerate hermitian forms 
$h: M\x M \to A$, where $M$ is a finitely generated right $A$-module (cf. \cite[Chap.~I]{Knus} or 
\cite[Chap.~7]{Sch}). We identify hermitian forms with their
Witt class in $W(A,\s)$, unless indicated otherwise.
Given an ordering $P\in X_F$ we wish
to define a signature at $P$, i.e. a morphism of groups
\[W(A,\s) \to \Z.\] 
Following the approach of \cite{BP2} we do this by extending scalars to a real closure $F_P$ of $F$ at $P$ and realizing 
that, by Morita equivalence, 
the Witt group of any $F_P$-algebra with involution is isomorphic to either $\Z$,  $0$ or $\Z/2\Z$. 
In the last two cases, the only sensible definition is to take the signature at $P$ to be identically zero. In this
case we call $P$
a \emph{nil-ordering}  and we write $\Nil[A,\s]$ for the set of all nil-orderings, noting that it only depends on
the Brauer class of $A$ and the type of $\s$. Furthermore, $\Nil[A,\s]$ is clopen in $X_F$, cf. 
\cite[Corollary~6.5]{A-U-Kneb}.

In the first case, the Witt group $W(A\ox_F F_P, \s\ox\id)$ is isomorphic to one of $W(F_P)$, $W((-1,-1)_{F_P}, -)$ 
or $W(F_P(\sqrt{-1}), -)$, where $-$ denotes (quaternion) conjugation, each one in turn being isomorphic to $\Z$ via
the usual Sylvester signature of quadratic or hermitian forms. The composite map $s_P$, given by
\begin{equation*}
\xymatrix{W(A,\s) \ar[r]& W(A\ox_F F_P, \s\ox\id)\ar[r] & \Z
},
\end{equation*}
enables us to define a signature. The map $s_P$ is independent of the choice of $F_P$ \cite[Prop.~3.3]{A-U-Kneb},  but a different
choice of Morita equivalence may result in a sign change  \cite[Prop.~3.4]{A-U-Kneb} and, conversely, such a sign change can always be obtained
by taking a well-chosen different Morita equivalence. 

At first sight, one way to fix a sign would be to demand that $s_P(\qf{1}_\s)$ is positive, as is the case for quadratic forms.
This is the approach
taken in \cite{BP2}, but it  may not always work,
since it may happen that $s_P(\qf{1}_\s)$ is in fact $0$, as illustrated in \cite[Rem.~3.11 and Ex.~3.12]{A-U-Kneb}.
Our solution to this dilemma is to show that there exists a hermitian form $\eta$ over $(A,\s)$, called a
\emph{reference form}, such that $s_P(\eta)$  is always nonzero whenever $P\in \wt X_F:=X_F\sm \Nil[A,\s]$, cf.
\cite[Prop.~3.2]{A-U-prime}.
Using this, given $P\in \wt X_F$, we define  the \emph{signature at $P$ with respect to the reference form $\eta$},
\[\sign_P^\eta: W(A,\s) \to \Z,\] 
to be the  map 
$s_P$, multiplied by $-1$ in case $s_P(\eta)<0$, so that $\sign_P^\eta (\eta)>0$.

The map $\sign_P^\eta$ does not depend on the Morita equivalence used in its computation
and so we may use the
explicit Morita equivalence presented in \cite{LU2} in all practical situations.

\begin{remark} In case $(A,\s)= (F, \id_F)$, we may take $\eta=\qf{1}$ and $\sign_P^\eta$ is then the usual
Sylvester signature $\sign_P$ of quadratic forms.
\end{remark}

\begin{remark} The signature map is defined for all hermitian forms over $(A,\s)$, not just the nondegenerate ones
as the notation above (which makes use of  $W(A,\s)$) might suggest. 
It suffices to replace a form by its nondegenerate part (cf. \cite[Prop.~A.3]{A-U-PS}), or alternatively, to replace $W(A,\s)$ by $\mathfrak{Herm}(A,\s)$,
the category of hermitian forms over $(A,\s)$.

\end{remark}

\begin{remark} A reference \emph{tuple} of hermitian forms of dimension one can be used instead of the reference form,
 cf.  \cite[Thm.~6.4]{A-U-Kneb} and \cite[\S 3]{A-U-prime}. In fact, this is the approach used in \cite{A-U-Kneb}.
\end{remark}

We collect some immediate properties of the signature map: 

\begin{prop}[Properties of the signature map {\cite[Thm~2.6]{A-U-prime}}]\mbox{}\label{prop2.4}
\begin{enumerate}[\ \ \ \ $(1)$]
\item Let $h$ be a hyperbolic form over $(A,\s)$, then
$\sign^\eta_P h=0$.

\item Let $h_1, h_2 \in W(A,\s)$, then
$\sign^\eta_P (h_1\perp h_2)=\sign^\eta_P h_1 + \sign^\eta_P h_2$.

\item Let $h\in W(A,\s)$  and $q\in W(F)$, then
$\sign_P^\eta (q\cdot h) = \sign_P q \cdot \sign_P^\eta h$.

\item \textup{(}Going-up\textup{)} Let $h\in W(A,\s)$ and let $L/F$ be an algebraic extension of ordered fields. Then 
\[\sign_Q^{\eta\ox L} (h\ox L)=\sign_{Q\cap F}^\eta h\] 
for all $Q\in X_L$. 

\end{enumerate}
\end{prop}

Property (4) is complemented by the following \emph{going-down} result:

\begin{thm}[Knebusch trace formula {\cite[Thm~8.1]{A-U-Kneb}}] 
Let $L/F$ be a finite extension of ordered fields and assume $P\in X_F$ extends to $L$.
Let $h\in W(A\ox_F L,\break \s\ox\id)$. Then
\[\sign_P^\eta(\mathrm{Tr}^*_{A\ox_F L} h) = \sum_{P\subseteq Q \in X_L} \sign_Q^{\eta\ox L} h,\]
 where $\mathrm{Tr}^*_{A\ox_F L} h$ denotes the Scharlau transfer induced by the $A$-linear homomorphism $\id_A \ox \mathrm{Tr}_{L/F}: A\ox_F L\to A$.
\end{thm}

\begin{thm}[Preservation under Morita equivalence  {\cite[Thm~4.2]{A-U-prime}}]
Let $(B,\tau)$ be an $F$-algebra with involution, Morita equivalent to  $(A,\s)$, and assume that $\s$ and $\tau$ are of the same  type. Let $\zeta: W(A,\s)\simtoo W(B,\tau)$ be the induced isomorphism of Witt groups. Then
\[\sign_P^\eta h = \sign_P^{\zeta(\eta)} \zeta(h)\]
for all $h\in W(A,\s)$ and all $P\in X_F$.
\end{thm}

\begin{thm}[Pfister's local-global principle {\cite[Thm~4.1]{LU1}}] 
Let $h\in W(A,\s)$. Then $h$ is a torsion form if and only if $\sign^\eta_P h=0$ for all $P\in X_F$.
\end{thm}

\begin{thm}[Continuity of the total signature {\cite[Thm~7.2]{A-U-Kneb}}]\label{cont} 
Let $h\in W(A,\s)$. The total signature of $h$, $\sign^\eta_\bullet h:X_F\to \Z, P\mapsto \sign_P^\eta h$, is continuous.
\end{thm}

These two theorems motivate the following results, familiar from the quadratic forms case. Let 
$ C(X_F,\mathbb{Z})_{[A,\s]}$ denote  the ring of continuous functions from $X_F$ to $\mathbb{Z}$, that are zero 
on $\Nil[A,\s]$.

\begin{thm}[{\cite[Prop.~4.3]{A-U-stab}}]
For every $f\in C(X_F,\mathbb{Z})_{[A,\s]}$ there exists $n\in\mathbb{N}$ such that $2^n f\in \im \sign^\eta$. In other words, the cokernel of $\sign^\eta$
is a $2$-primary torsion group.
\end{thm}

 The \emph{stability index} of $(A, \sigma )$ is the smallest $k\in \mathbb{N}$ such that $2^k C(X_F,\mathbb{Z})_{[A,\s]} \subseteq \im \sign^\eta$ if such a $k$ exists and $\infty$ otherwise. It is independent of the choice of $\eta$.
The group $\coker \sign^\eta$ is up to isomorphism independent of the choice of $\eta$. We denote it by $S_\eta(A, \sigma )$ and call it
the \emph{stability group} of $(A, \sigma )$.

\begin{thm}[{\cite[Thm.~4.10]{A-U-stab}}] 
Let $W_t(A, \sigma )$ denote the torsion subgroup of $W(A, \sigma )$. \textup{(}Recall that it is $2$-primary torsion by \cite[Thm.~5.1]{Sch1}.\textup{)}
The sequence
\begin{equation*}
\xymatrix{
0 \ar[r] &   W_t(A, \sigma ) \ar[r] & W(A, \sigma) \ar[r]^--{\sign^\eta} &  C(X_F,\mathbb{Z})_{[A,\s]} \ar[r] & S_\eta(A, \sigma )\ar[r] & 0
}
\end{equation*}
is exact. The group $S_\eta(A, \sigma )$ is a  $2$-primary torsion group.
\end{thm}

\section{Ideals and morphisms}

Let $R$ be a commutative ring and let $M$ be an $R$-module. We introduce ideals of $R$-modules as follows:
An \emph{ideal} of $M$ is a pair $(I,N)$ where $I$ is an ideal of $R$ and $N$
      is a submodule of $M$ such that $I \cdot M \subseteq N$.
An ideal $(I,N)$ of $M$ is \emph{prime} if $I$ is a prime ideal of $R$ (we
      assume that all prime ideals are proper), $N$ is a proper submodule of
      $M$, and         
      for every $r \in R$ and $m \in M$, $r\cdot m \in N$ implies that $r \in I$
      or $m \in N$.

These definitions are in part motivated by the following natural example:
The pair $(\ker \sign_P, \ker \sign_P^\eta)$ is a prime ideal of the $W(F)$-module $W(A, \sigma )$
whenever $P \in \widetilde{X}_F$ .

We obtain a classification \`a la Harrison and Lorenz-Leicht \cite{Lorenz-Leicht}:

\begin{thm}[{\cite[Props.~6.5, 6.7]{A-U-prime}}]
\label{classif}
  Let $(I,N)$ be a prime ideal of the $W(F)$-module $W(A, \sigma )$. 
  \begin{enumerate}[\ \ \ \ $(1)$]
  \item If $2 \not \in I$, then one of
  the following holds:
  \begin{enumerate}[$(i)$]
    \item There exists   $P \in \wt X_F$ such that $(I,N) = (\ker \sign_P,  \ker \sign_P^\eta)$.
    \item There exist   $P \in \wt X_F$ and a prime $p >2$   such that
    $(I,N) = \bigl(\ker (\pi_p \circ \sign_P), \ker (\pi \circ \sign_P^\eta) \bigr)$,
     where $\pi_p: \mathbb{Z}\to \mathbb{Z}/p\mathbb{Z}$ and  
      $\pi: \im \sign_P^\eta \to \im \sign_P^\eta/(p\cdot \im\sign_P^\eta)$ are the canonical projections.
  \end{enumerate}
  
  \item If $2\in I$, then $I=I(F)$, the fundamental ideal of $W(F)$. Furthermore, 
  a pair $(I(F), N)$ is a prime ideal of $W(A, \sigma )$ if and only if $N$ is a proper submodule of $W(A, \sigma )$ with 
$I(F)\cdot W(A, \sigma )\subseteq N$.
    \end{enumerate}
\end{thm}

\begin{remark} When $2\not\in I$, $N$ is completely determined by $I$. This is however not the case when $2\in I$, cf. \cite[Ex.~6.8]{A-U-prime}.
\end{remark}

There is a notion of morphism linked in the usual way to the above notion of ideal, cf. 
\cite[Lemmas~5.6, 5.7 and 5.8]{A-U-prime}:
 Let $R$ and $S$ be commutative rings, let $M$ be an $R$-module and $N$  an
  $S$-module. We say that a pair $(\vf, \psi)$ is an \df{$(R,S)$-morphism 
  \textup{(}of modules\textup{)} from 
  $M$ to $N$}
  if 
  \begin{enumerate}[\ \ \ \ $(1)$]
    \item $\vf : R \rightarrow S$ is a morphism of rings (and in particular $\vf(1)=1$);
    \item $\psi : M \rightarrow N$ is a morphism of additive groups;
    \item for every $r \in R$ and $m \in M$, $\psi(r\cdot m) = \vf(r)\cdot \psi(m)$.
  \end{enumerate}
  We call an $(R,S)$-morphism $(\vf,\psi)$  \df{trivial} if $\psi=0$. 
  We denote  the set of all $(R,S)$-morphisms from $M$ to
  $N$ by $\Hom_{(R,S)}(M,N)$ and its subset of nontrivial $(R,S)$-morph\-isms by $\Hom^*_{(R,S)}(M,N)$.

  Let $(\vf, \psi_1)$ and $(\vf, \psi_2)$ be $(R,S)$-morphisms of modules from
  $M$ to $N_1$ and $N_2$ respectively. We say that $(\vf, \psi_1)$ and $(\vf,
  \psi_2)$ are \df{equivalent} if there is an isomorphism of $\im \vf$-modules
  $\vartheta : \im \psi_1
  \rightarrow \im \psi_2$ such that $\psi_2 = \vartheta \circ \psi_1$. We write
  $\sim$ for the relation ``being equivalent''.

The pair
$(\sign_P, \sign_P^\eta)$ is again a natural example of a $(W(F),\Z)$-morphism from
    $W(A,\s)$ to $\Z$ and is trivial if and only if $P\in \Nil[A,\s]$.

The classification of prime ideals of $W(A,\s)$ yields the following description of signatures as morphisms:

\begin{thm}[{\cite[Prop.~7.4]{A-U-prime}}]\label{morph}
  The map that sends $P \in \wt{X}_F$ to the pair $(\sign_P,\break \sign^\eta_P)$ induces a
  bijection between $\wt{X}_F$ and the equivalence classes with respect to $\sim$ of
  $\Hom^*_{(W(F),\Z)}(W(A,\s),\Z)$.
\end{thm}

Theorems~\ref{classif} and \ref{morph} give us

\begin{cor} There is a bijective correspondence between the prime ideals of $W(A,\s)$, the nonzero signatures
of hermitian forms over $(A,\s)$ and  the equivalence classes of  $\Hom^*_{(W(F),\Z)}(W(A,\s),\Z)$ 
with respect to $\sim$.
\end{cor}

  \section{Sums of hermitian squares}\label{shs}
In the field case, Pfister's local-global principle can be used to give a short
proof of the fact that sums of squares are exactly the elements that are
nonnegative at every ordering. In \cite{A-U-PS} we showed that the same approach directly yields a similar
result for $F$-division algebras with involution and, with some extra effort, for
all $F$-algebras with involution.

Let $A^\x$ denote the set of invertible elements of $A$,
$\Sym(A,\s)$ the set of $\s$-symmetric elements of $A$ and $\sas^\x:=\sas \cap A^\x$. 
We say that an element $a \in \Sym(A,\s)$ is \emph{$\eta$-maximal at an ordering $P
\in X_F$} if $\sign^\eta_P \qf{a}_\s$ is maximal among all $\sign^\eta_P
\qf{b}_\s$ for $b \in \Sym(A,\s)$. In the field case, this means $\sign_P
\qf{a} = 1$, in other words $a \in P \setminus \{0\}$. For elements  $b_1,\ldots, b_t \in F^\x$ we 
denote the \emph{Harrison set} $\{P \in X_F \mid b_1,\ldots, b_t \in P\}$ by $H(b_1,\ldots, b_t)$.

\begin{thm}[{\cite[Thm.~3.6]{A-U-PS}}]
\label{main_thm_2}
  Let  $b_1,\ldots, b_t \in F^\x$ and $Y=H(b_1,\ldots, b_t)$. Assume that
  $a\in\sas^\x$ is $\eta$-maximal at all  $P \in Y$.  Let $u\in\Sym(A,\s)$.
The
  following statements are equivalent:
  \begin{enumerate}[\ \ \ \ $(i)$]
    \item $u$ is $\eta$-maximal  at all $P\in Y$.
    \item $u \in \das (k \x \pf{b_1,\ldots, b_t}  \qf{a}_\s)$ for some $k \in \N$.
  \end{enumerate}
\end{thm}
The presence of the element $a$ as well as the hypothesis on $\eta$-maximality correspond in the field case to the
fact  that $1$ belongs to every ordering. Here $1$ does not
play a particular role since it may not have maximal signature at some
orderings. We replace it by the element $a$ and only consider a set of
orderings $Y$ on which $a$ has maximal signature.

\begin{remark}\label{max}  
As a consequence of our study of positive involutions,
given $Q\in \wt X_F$, there always exist $a$ and $Y$ that satisfy the hypothesis of 
Theorem~\ref{main_thm_2} with $Q\in Y$, cf. Theorem~\ref{positive=max} below, together with 
Theorem~\ref{cont} on the continuity of the signature map.
\end{remark}

Consider the form $\tasu(x,y):= \Trd_A (\s(x)uy)$ for $x,y\in A$ and,
following \cite[Def.~1.1]{P-S}, let
\begin{align*}
  X_\s &:= \{P \in X_F \mid \s \text{ is positive at } P\} \\
       &\phantom{:}= \{P \in X_F \mid T_{(A,\s,1)} \textrm{ is positive semidefinite at } P\}.
\end{align*}
Theorem~\ref{main_thm_2}
is reminiscent of Procesi and Schacher's theorem \cite[Thm.~5.4]{P-S}
 and, with the notation just introduced, can
also be used to address the question they raised in \cite[p.~404]{P-S}, 
of whether or not the following property is always true:
\begin{description}
\setlength{\labelwidth}{2em}
  \item[(PS)] for every $u\in \Sym(A,\s)$, the form $\tasu$ is positive
    semidefinite at all  $P \in X_\s$ if and only if $u \in \das (s \x
    \qf{1}_\s)$ for some $s \in \N$.
\end{description}

The general answer to this question is negative as shown in
\cite{K-U}, but we can now describe cases where the answer is
positive, and also propose a natural reformulation (inspired by signatures of hermitian forms)
of the question that has a positive answer.

\begin{prop}[{\cite[Cor.~4.18]{A-U-PS}}]
\label{canary}
  If $X_\s =\wt X_F$, then property \textup{(PS)} holds.
\end{prop}
And, if we introduce the property
\begin{description}
\setlength{\labelwidth}{2em}
  \item[(PS')] for every $u\in \Sym(A,\s)$, the form $\tasu$ is positive
    semidefinite at all  $P \in \wt X_F$ if and only if $u \in \das (s \x
    \qf{1}_\s)$ for some $s \in \N$,
\end{description}
we obtain
\begin{thm}[{\cite[Thm.~4.19]{A-U-PS}}]
\label{biscuit}
  Property \textup{(PS')} holds if and only
if $\wt X_F = X_\s$.
\end{thm}

\section{Positive cones}
The results presented thus far suggest that there could be a notion of
``ordering'' on central simple algebras with involution, whose behaviour would
be similar to that of orderings on fields. The purpose of this final section is
to present such a notion.

\begin{defi}[{\cite[Def.~3.1]{A-U-pos}}]
\label{def-preordering}
  A \emph{prepositive cone} $\CP$ on $(A,\s)$ is
  a subset $\CP$ of $\Sym(A,\s)$  such that
  \begin{enumerate}[\ \ \ \ (P1)]
    \item $\CP \not = \varnothing$;
    \item $\CP + \CP \subseteq \CP$;
    \item $\s(a) \cdot \CP \cdot a \subseteq \CP$ for every $a \in A$;
    \item $\CP_F := \{u \in F \mid u\CP \subseteq \CP\}$ is an ordering on $F$.
    \item $\CP \cap -\CP = \{0\}$ (we say that $\CP$ is \emph{proper}).
  \end{enumerate}
  We say that a prepositive cone $\CP$ is \emph{over $P\in X_F$}
if $\CP_F=P$. 
\end{defi}
\begin{remark}
  Axiom (P4) is necessary if we want our prepositive cones to consist of either
  positive semidefinite (PSD) matrices with respect to $P$, or of 
  negative semidefinite (NSD) matrices with respect to $P$,
 in the case of $(M_n(F),t)$, see \cite[Rem.~3.13]{A-U-pos}.
 
  If $\CP$ is a prepositive cone, then $-\CP$ is also a prepositive cone. This is
  due to the fact that prepositive cones are meant to contain elements of
  maximal signature, and the sign of the signature can vary with a change of the
  reference form.
\end{remark}

It can be shown that there is a prepositive cone over $P \in X_F$ on $(A,\s)$
if and only if $P \in \wt X_F$, cf. \cite[Prop.~6.6]{A-U-pos}.

\begin{examples}\mbox{}
  \begin{enumerate}[\ \ \ \ $(1)$]
    \item  Let $P\in \wt X_F$. We define
      \[\CM_P^\eta(A,\s):=\{ a \in \Sym(A,\s)^\x \mid a \text{ is $\eta$-maximal
      at } P\} \cup \{0\}.\]
      If $A$ is a division algebra, $\CM_P^\eta(A,\s)$ is a prepositive cone over
      $P$ on $(A,\s)$.
    \item The set of PSD matrices, and the set  of NSD matrices with respect to some $P
      \in X_F$ are both prepositive cones over $P$ on $(M_n(F),t)$.
  \end{enumerate}
 
\end{examples}

\begin{remark}
Other notions of orderings have been introduced for 
division rings with involution, most notably Baer orderings, $*$-orderings and their variants and an extensive
theory has been
developed around them. Craven's surveys \cite{Cr1} and \cite{Cr2}  provide more information on these topics. 
Without going into the details, the main difference in the definitions  is that positive cones were developed to correspond 
to a pre-existing algebraic notion, namely signatures of hermitian forms (e.g. axiom (P4) reflects the fact that the 
signature is a morphism of modules, cf.  Proposition~\ref{prop2.4}(3); see also the sentence after Theorem~\ref{thm5.8})
and as a consequence are not required to induce total orderings on the set of symmetric elements.
\end{remark}

We obtain the desired results linking prepositive cones and
$W(A,\s)$:

\begin{prop}[{\cite[Prop.~7.11]{A-U-pos}}]
\label{tor_fr}
 The following statements are equivalent:
  \begin{enumerate}[\ \ \ \ $(i)$]
    \item  $(A,\s)$ is formally real \textup{(}i.e. there is at least one prepositive cone
      on $(A,\s)$\textup{)};
    \item $W(A,\s)$ is not torsion;
    \item $\wt X_F\not=\varnothing$.
 \end{enumerate}
\end{prop}

Prepositive cones are well-behaved under Morita equivalence: If $(A,\s)$ and
$(B,\tau)$ are Morita equivalent, then there is an inclusion-preserving bijection between their sets
of prepositive cones. This bijection can be made explicit in the case of scaling,
or when $(A,\s) = (D,\vt)$ is an $F$-division algebra with involution  
and $(B,\tau) = (M_n(D),\vt^t)$, using descriptions
of prepositive cones reminiscent of the characterizations of positive semidefinite
matrices:

\begin{prop}[{\cite[\S4.1, \S4.2]{A-U-pos}}]\label{explicit}
The prepositive cones on $(M_n(D),\vt^t)$ are of the form 
\[\PSD_n(\CP): = \{M \in \Sym(M_n(D), \vt^t) \mid \forall X \in D^n \quad \vt(X)^t M
X \in \CP\},\]
where $\CP$ is a prepositive cone on $(D,\vt)$. 

The prepositive cones on $(D,\vt)$ are of the form 
\[\Tr_n(\CP) :=  \{\tr(M) \mid M \in \CP\},\]
where $\CP$ is a prepositive cone on $(M_n(D),\vt^t)$. 
\end{prop}

We use prepositive cones to consider the question of the existence  of positive involutions:

\begin{thm}[{\cite[Thm.~6.8]{A-U-pos}}]\label{main_pos} 
Let  $P \in X_F$. The following statements are equivalent:
\begin{enumerate}[\ \ \ \ $(i)$]
\item There is an involution  $\tau$ on $A$ which is positive at $P$ and of the same type as $\s$;
\item $P \in \wt X_F = X_F\setminus \Nil[A,\s]$.
\end{enumerate}
\end{thm}

The notion of prepositive cone can be seen as somewhat equivalent to that of
preordering or Prestel's pre-semiordering \cite{Prestel73}, \cite{Prestel84}, so it is natural to consider in more detail the
maximal prepositive cones, which we simply call \emph{positive cones}. 
They can be completely described and match the examples
provided above. To see this we define, for $P \in X_F$ and $S \subseteq \sas$,
\[\CC_P(S) := \Bigl\{\sum_{i=1}^k u_i \s(x_i)s_ix_i \,\Big|\, k \in \N, \ u_i
\in P,\ x_i
  \in A,\ s_i \in S\Bigr\}\]
(the smallest, possibly nonproper, prepositive cone over $P$ containing $S$),
and we denote by $X_{(A,\s)}$ the set of all positive cones on $(A,\s)$.

\begin{thm}[{\cite[Thm.~7.5]{A-U-pos}}]\label{thm5.8}
\label{positive=max}

  \[X_{(A,\s)} =
  \{-\CC_P(\CM^\eta_P(A,\s)), \CC_P(\CM^\eta_P(A,\s)) \mid P \in  \wt X_F  \}.\]

  Moreover, for each $\CP \in X_{(A,\s)}$, there exists $\ve\in\{-1,1\}$ such that
  $\CP \cap A^\x = \ve \CM^\eta_P(A,\s) \setminus \{0\}$.
\end{thm}

In particular, the only positive cones over $P$ on $(D,\vt)$ are 
$\CM^\eta_P(D,\vt)$ and $-\CM^\eta_P(D,\vt)$ and therefore the examples above 
are essentially the only positive cones on $(A,\s)$, cf. \cite[Props.~4.3 and 4.9]{A-U-pos}. 
It follows that the PSD matrices over $P$
and the NSD matrices over $P$ are the only positive cones over $P$ on
$(M_n(F),t)$. (See also Proposition~\ref{explicit}.)

Using this description, it is possible to make the link with the results
presented in Section~\ref{shs}, and to obtain results similar to the
Artin-Schreier and Artin theorems.

\begin{thm}[{\cite[Thm.~7.9]{A-U-pos}}]
\label{fr}
  The following statements are equivalent:

  \begin{enumerate}[\ \ \ \ $(i)$]
    \item\label{fr1} $(A,\s)$ is formally real;

    \item\label{fr3} There is $a \in \Sym(A,\s)^\x$ and $P \in X_F$ such that 
      $\CC_P(a) \cap -\CC_P(a)=\{0\}$;
    \item\label{fr4} There is $b \in \Sym(A,\s)^\x$ such that $\qf{b}_\s$ is 
      strongly anisotropic \textup{(}i.e. the forms $k\x \qf{b}_\s$ have no  nontrivial zeroes for all $k\in\N$\textup{)}.

  \end{enumerate}
\end{thm}
The second statement is a trivial consequence of the first one, but it is
still included here to point out that while the element $a$ in it obviously
belongs to a prepositive cone (namely $\CC_P(a)$), the element $b$ in the third
statement may not belong to any prepositive cone on $(A,\s)$, contrary to what
could be expected from the field case (see {\cite[Rem.~7.10]{A-U-pos}}).

\begin{thm}[{\cite[Thm.~7.14]{A-U-pos}}]
\label{intersection}
  Let $b_1,\ldots,b_t \in F^\x$,
  let $Y =
  H(b_1,\ldots,b_t)$ and let $a \in \Sym(A,\s)^\x$ be such that, for
  every
  $\CP \in X_{(A,\s)}$ with   $\CP_F \in Y$, $a \in \CP\cup -\CP$. Then
  \[\bigcap \{\CP \in X_{(A,\s)} \mid \CP_F \in Y \textrm{ and } a \in \CP\} =
  \bigcup_{s \in \N} D_{(A,\s)} (s \x \pf{b_1,\ldots,b_t}
  \qf{a}_\s).\]
\end{thm}

As a consequence of our study of positive involutions,
given $Q\in \wt X_F$, there always exist $a$ and $Y$ that satisfy the hypothesis of 
Theorem~\ref{intersection} with $Q\in Y$, cf. Remark~\ref{max}.

The element $a$ in this theorem plays the same role as the element $a$ in
Theorem~\ref{main_thm_2}, and chooses a prepositive cone from $\{\CP, -\CP\}$
in a uniform way. This is not necessary in the field case, because $1$ belongs to
every ordering.
 In the
special case where  $a=1$ can be used for this purpose, we obtain a
result more similar to the usual one:
\begin{cor}[{\cite[Cor.~7.15]{A-U-pos}}]
\label{Artin}
Assume that for every $\CP \in X_{(A,\s)}$, $1\in
  \CP\cup -\CP$. Then
  \[\bigcap \{\CP \in X_{(A,\s)} \mid 1\in \CP\} =
  \Bigl\{\sum_{i=1}^s \s(x_i)x_i \,\Big|\, s\in \N,\ x_i \in
  A\Bigr\}.\]
\end{cor}
The hypothesis of Corollary~\ref{Artin} is exactly $X_\s = \wt X_F$ in the
terminology of Section~\ref{shs}. More precisely, as seen therein, this property characterizes the
algebras with involution for which there is a positive answer 
to (PS'),  cf. \cite[Section~4.2]{A-U-PS}.

A natural question is to ask if signatures of hermitian forms over $(A,\s)$ can now also be defined with respect to
positive cones on $(A,\s)$. As shown in \cite[\S 8.2]{A-U-pos}, this can indeed be done using decompositions of hermitian
forms, reminiscent of Sylvester's decomposition for quadratic forms:

\begin{thm}[{\cite[Cor.~8.14, Lemma~8.15]{A-U-pos}}]  
There exists an integer $t$, depending on $(A,\s)$, such that for every $\CP \in X_{(A,\s)}$ and for every
$h\in W(A,\s)$ there exist $u_1,\ldots, u_t \in P:=\CP_F $, $a_1,\ldots, a_r \in \CP \cap A^\x$ and $b_1,\ldots, b_s
\in -\CP\cap A^\x$ such that
\[n_P^2\x \qf{u_1,\ldots, u_t}\ox h \simeq \qf{a_1,\ldots, a_r}_\s \perp \qf{b_1,\ldots, b_s}_\s,\]
where $n_P$ is the matrix size of $A\ox_F F_P$, and
$r$ and $s$ are positive integers, uniquely determined by $\CP$ and the rank of $h$.
\end{thm}

This theorem allows us to define
\[\sign_\CP h:= \frac{r-s}{n_Pt} \in \Z,\]
cf. \cite[Def.~8.16]{A-U-pos}, which
coincides with the signature defined in Section~\ref{sec:two}, cf. \cite[Prop.~8.17]{A-U-pos} and also yields 
total signature maps that are continuous for the topology defined below, cf. \cite[Thm.~9.19]{A-U-pos}.

We finish this paper by a presentation of our main results concerning the
topology of $X_{(A,\s)}$. We define, for $a_1, \ldots, a_k \in
\Sym(A,\s)$,
\[H_\s(a_1, \ldots, a_k) := \{\CP \in X_{(A,\s)} \mid a_1, \ldots,
a_k \in
\CP\}.\]
We denote by $\CT_\s$ the topology on $X_{(A,\s)}$ generated by
the sets
$H_\s(a_1, \ldots, a_k)$, for $a_1, \ldots, a_k \in \Sym(A,\s)$, and
by
$\CT_\s^\x$ the topology on $X_{(A,\s)}$ generated by the sets
$H_\s(a_1,
\ldots, a_k)$, for $a_1, \ldots, a_k \in \Sym(A,\s)^\x$.

\begin{prop}[{\cite[Prop.~9.2]{A-U-pos}}]
\label{top_equal}
  The topologies $\CT_\s$ and $\CT_\s^\x$  are equal.
\end{prop}

Recall that spectral topologies, defined in \cite{Hochster}, are precisely the topologies of the spectra of commutative 
rings,
and that a map between spectral spaces (i.e. spaces equipped with spectral topologies) is called spectral if it is 
continuous and the preimage of a quasicompact open set is quasicompact.

\begin{thm}[{\cite[Prop.~9.17]{A-U-pos}}]
\label{spectral-top}
  $\CT_\s$ is a spectral topology on $X_{(A,\s)}$.
\end{thm}

The topology $\CT_\s$ is also well-behaved under Morita equivalence:
\begin{prop}[{\cite[Prop.~9.18]{A-U-pos}}]
\label{MMMMorita}
 Let $(A,\s)$ and $(B,\tau)$ be two Morita
  equivalent 
  $F$-algebras with involution. The spaces $(X_{(A,\s)}, \CT_\s)$
  and $(X_{(B,\tau)}, \CT_\tau)$ are homeomorphic via a spectral map.
\end{prop}

\def\cprime{$'$}


\begin{thebibliography}{33}

\bibitem{Albert35}
  \newblock A.A. Albert, 
  \newblock {Involutorial simple algebras and real {R}iemann matrices},
  \newblock  \emph{Ann. of Math. (2)} \textbf{36} (1935), no.~4, 886--964.

\bibitem{Artin}
  \newblock E.~Artin, 
  \newblock {\"{U}ber die {Z}erlegung definiter {F}unktionen in {Q}uadrate},
  \newblock \emph{Abh. Math. Sem. Univ. Hamburg} \textbf{5} (1927), no.~1, 100--115.

\bibitem{Artin-Schreier} 
  \newblock E.~Artin and O.~Schreier, 
  \newblock {Algebraische {K}onstruktion reeller {K}\"orper}, 
  \newblock \emph{Abh. Math. Sem. Univ. Hamburg} \textbf{5} (1927), no.~1, 85--99.

\bibitem{A-U-Kneb} 
  \newblock V.~Astier and T.~Unger, 
  \newblock {Signatures of hermitian forms and the {K}nebusch  trace formula}, 
  \newblock \emph{Math. Ann.} \textbf{358} (2014), no.~3-4, 925--947.

\bibitem{A-U-prime} 
  \newblock V.~Astier and T.~Unger, 
  \newblock {Signatures of hermitian forms and ``prime ideals'' of {W}itt groups}, 
  \newblock \emph{Adv. Math.} \textbf{285} (2015), 497--514.

\bibitem{A-U-pos}
  \newblock V.~Astier and T.~Unger, 
  \newblock {Positive cones on algebras with involution},
  \newblock preprint (2017), \url{http://arxiv.org/abs/1609.06601}.

\bibitem{A-U-PS}
  \newblock V.~Astier and T.~Unger, 
  \newblock {Signatures of hermitian forms, positivity, and an answer to a
  question of {P}rocesi and {S}chacher}, 
  \newblock preprint (2016),  \url{http://arxiv.org/abs/1511.06330}.


\bibitem{A-U-survey2017}
  \newblock V.~Astier  and T.~Unger, 
  \newblock {Signatures, sums of hermitian squares and positive cones on algebras with involution},
  \newblock \emph{S{\'e}minaire de Structures Alg{\'e}briques Ordonn{\'e}es 2015--2016} \textbf{91} (2017).                                                                                                                                                                                                                                                                                                                                                                                                                        


\bibitem{A-U-stab} 
  \newblock V.~Astier and T.~Unger, 
  \newblock {Stability index of algebras with involution}, 
  \newblock \emph{Contemporary Mathematics} \textbf{697} (2017), 41--50.
  

\bibitem{auel-2011} 
  \newblock A.~Auel, E.~Brussel, S.~Garibaldi, and U.~Vishne, 
  \newblock {Open problems on central simple algebras}, 
  \newblock \emph{Transform. Groups} \textbf{16} (2011), no.~1, 219--264.

\bibitem{BP2} 
  \newblock E.~Bayer-Fluckiger and R.~Parimala, 
  \newblock {Classical groups and the {H}asse principle}, 
  \newblock \emph{Ann. of Math. (2)} \textbf{147} (1998), no.~3, 651--693.

\bibitem{BCR} 
  \newblock J.~Bochnak, M.~Coste, and M.-F. Roy, 
  \newblock \emph{Real algebraic geometry}, 
  \newblock Ergebnisse der Mathematik und ihrer Grenzgebiete (3), vol.~36, Springer-Verlag, Berlin,
  1998.

\bibitem{Cr1} 
  \newblock T.C. Craven,
  \newblock {Orderings, valuations, and Hermitian forms over $*$-fields}, 
  \newblock in \emph{K-theory and algebraic geometry: connections with quadratic forms and division algebras (Santa Barbara, CA, 1992)}, 
  \newblock {Proc. Sympos. Pure Math.} {58}, Part 2, American Mathematical Society, Providence, RI, (1995), 149--160.  

\bibitem{Cr2} 
  \newblock T.C. Craven, 
  \newblock {Valuations and Hermitian forms on skew fields}, 
  \newblock in \emph{Valuation theory and its applications, Vol. I (Saskatoon, SK, 1999)},  
  \newblock Fields Inst. Commun., 32, American Mathematical Society, Providence, RI, (2002), 103--115.

\bibitem{Hilbert1900}
  \newblock D.~Hilbert, 
  \newblock {Mathematische Probleme. Vortrag, gehalten auf dem internationalen Mathe\-matiker-Congress zu Paris 1900.}, 
  \newblock \emph{Nachr. Ges. Wiss. G\"ottingen, Math.-Phys. Kl.}
  \textbf{1900} (1900), 253--297.

\bibitem{Hochster} 
  \newblock M.~Hochster, 
  \newblock {Prime ideal structure in commutative rings}, 
  \newblock \emph{Trans. Amer.  Math. Soc.} \textbf{142} (1969), 43--60.

\bibitem{K-U} 
  \newblock I.~Klep and T.~Unger, 
  \newblock {The {P}rocesi-{S}chacher conjecture and {H}ilbert's 17th problem for algebras with involution}, 
  \newblock \emph{J. Algebra} \textbf{324} (2010), no.~2, 256--268.

\bibitem{Knus}
  \newblock M.-A. Knus, 
  \newblock \emph{Quadratic and {H}ermitian forms over rings}, 
  \newblock Grundlehren der Mathematischen Wissenschaften, vol. 294, Springer-Verlag, Berlin, 1991.

\bibitem{BOI} 
  \newblock M.-A. Knus, A.~Merkurjev, M.~Rost, and J.-P. Tignol, 
  \newblock \emph{The book of involutions}, 
  \newblock American Mathematical Society Colloquium Publications, vol.~44, American Mathematical Society, Providence, RI, 1998.

\bibitem{Lam-1984} 
  \newblock T.Y. Lam, 
  \newblock {An introduction to real algebra}, 
  \newblock in \emph{Ordered fields and real algebraic geometry (Boulder, Colo., 1983)},
  \newblock Rocky Mountain J. Math. \textbf{14} (1984), no.~4, 767--814, .

\bibitem{Lam} 
  \newblock T.Y. Lam, 
  \newblock \emph{Introduction to quadratic forms over fields}, 
  \newblock Graduate Studies in Mathematics, vol.~67, American Mathematical Society, Providence, RI, 2005.

\bibitem{LT} 
  \newblock D.W. Lewis and J.-P. Tignol, 
  \newblock {On the signature of an involution}, 
  \newblock \emph{Arch. Math. (Basel)} \textbf{60} (1993), no.~2, 128--135.

\bibitem{LU1} 
  \newblock D.W. Lewis and T.~Unger, 
  \newblock {A local-global principle for algebras with involution and {H}ermitian forms}, 
  \newblock \emph{Math. Z.} \textbf{244} (2003), no.~3, 469--477.

\bibitem{LU2} 
  \newblock D.W. Lewis and T.~Unger, 
  \newblock {Hermitian {M}orita theory: a matrix approach}, 
  \newblock \emph{Irish Math. Soc. Bull.} \textbf{62} (2008),  37--41.

\bibitem{Lorenz-Leicht} 
  \newblock F.~Lorenz and J.~Leicht, 
  \newblock {Die {P}rimideale des {W}ittschen {R}inges},
  \newblock \emph{Invent. Math.} \textbf{10} (1970), 82--88.

\bibitem{Pfister66} 
  \newblock A.~Pfister, 
  \newblock {Quadratische {F}ormen in beliebigen {K}\"orpern}, 
  \newblock \emph{Invent. Math.} \textbf{1} (1966), 116--132.

\bibitem{Prestel73} 
  \newblock A. Prestel,
  \newblock {Quadratische {S}emi-{O}rdnungen und quadratische {F}ormen.}
  \newblock \emph{Math. Z.} \textbf{133} (1973), 319--342.

\bibitem{Prestel84} 
  \newblock A. Prestel,
  \newblock \emph{Lectures on formally real fields}, 
  \newblock Lecture Notes in Mathematics, vol.~1093, Springer-Verlag, Berlin, 1984.


\bibitem{P-S}
  \newblock C.~Procesi and M.~Schacher, 
  \newblock {A non-commutative real {N}ullstellensatz and {H}ilbert's 17th problem}, 
  \newblock \emph{Ann. of Math. (2)} \textbf{104} (1976), no.~3, 395--406.

\bibitem{Q} 
  \newblock A.~Qu{\'e}guiner, 
  \newblock {Signature des involutions de deuxi\`eme esp\`ece},
  \newblock \emph{Arch. Math. (Basel)} \textbf{65} (1995), no.~5, 408--412.

\bibitem{Sch1} 
  \newblock W.~Scharlau,
  \newblock {Induction theorems and the structure of the {W}itt group},
  \newblock \emph{Invent. Math.} \textbf{11} (1970), 37--44.

\bibitem{Sch} 
  \newblock W.~Scharlau, 
  \newblock \emph{Quadratic and {H}ermitian forms}, 
  \newblock Grundlehren der Mathematischen Wissenschaften, vol. 270, Springer-Verlag, Berlin, 1985.

\bibitem{Sylvester1852}
  \newblock J.J. Sylvester, 
  \newblock {A demonstration of the theorem that every homogeneous
  quadratic polynomial is reducible by real orthogonal substitutions to the
  form of a sum of positive and negative squares}, 
  \newblock \emph{Philosophical Magazine}, \textbf{IV} (1852), 138--142.

\bibitem{Tignol-1998} 
  \newblock J.-P. Tignol, 
  \newblock {Algebras with involution and classical groups}, 
  \newblock in \emph{European {C}ongress of {M}athematics, {V}ol.\ {II} ({B}udapest, 1996)}, 
  \newblock Progr. Math., vol. 169, Birkh\"auser, Basel, (1998), 244--258.

\bibitem{Weil} 
  \newblock A.~Weil, 
  \newblock {Algebras with involutions and the classical groups}, 
  \newblock \emph{J. Indian Math. Soc. (N.S.)} \textbf{24} (1961), 589--623.

\end{thebibliography}
\end{document}